\documentclass[12pt]{amsart}
\usepackage{latexsym}
\usepackage[all]{xy}
\usepackage{amsfonts}
\usepackage{amsthm}
\usepackage{amsmath}
\usepackage{amssymb}
\newcommand{\cC}{{\mathcal C}}

\def\sw#1{{\sb{(#1)}}}

\def\<{{\langle}}
\def\>{{\rangle}}

\def\eps{\epsilon}

\def\note#1{{}}

\def\note#1{}
\def\M{{\bf M}}

\def\cC{{\mathcal C}}

\def\hom#1#2#3#4{{{}\sb{#1}{\rm Hom}\sb{#2}(#3,#4)}}
\def\lhom#1#2#3{{{}\sb{#1}{\rm Hom}(#2,#3)}}
\def\rhom#1#2#3{{{\rm Hom}\sb{#1}(#2,#3)}}

\def\rend#1#2{{{\rm End}\sb{#1}(#2)}}
\def\lrend#1#2#3{{{}\sb{#1}{\rm End}\sb{#2}(#3)}}

\def\beq{\begin{equation}}
\def\eeq{\end{equation}}
\def\DC{{\Delta_\cC}}
\def \eC{{\eps_\cC}}

\def\id{{I}}
\def\im{{\rm Im}}

\newcommand{\scr}[1]{\mbox{\scriptsize $#1$}}
\def\ut{\scr{\otimes}}
\def\ot{{\otimes}}

\def\End{\mbox{\rm End}\,} 
\def\rM{\varrho^{M}}

\newcounter{zlist}
\newenvironment{zlist}{\begin{list}{(\arabic{zlist})}{
\usecounter{zlist}\leftmargin2.5em\labelwidth2em\labelsep0.5em
\topsep0.6ex
\parsep0.3ex plus0.2ex minus0.1ex}}{\end{list}}

\newtheorem{thm}{Theorem}[section]
\newtheorem{conjecture}[thm]{Conjecture}

\def\Label#1{\label{#1}\ifmmode\llap{[#1] }\else 
\marginpar{\smash{\hbox{\tiny [#1]}}}\fi}
\def\Label{\label}

\newtheorem{proposition}{Proposition}[section]

\newtheorem{corollary}[proposition]{Corollary}
\newtheorem{theorem}[proposition]{Theorem}

\theoremstyle{definition}
\newtheorem{definition}[proposition]{Definition}
\newtheorem{example}[proposition]{Example}

\theoremstyle{remark}
\newtheorem{remark}[proposition]{Remark}

\newcounter{c}
\renewcommand{\[}{\setcounter{c}{1}$$}
\newcommand{\etyk}[1]{\vspace{-7.4mm}$$\begin{equation}\Label{#1}
\addtocounter{c}{1}}
\renewcommand{\]}{\ifnum \value{c}=1 $$\else \end{equation}\fi}
\begin{document}

\title{On coseparable and biseparable corings}
\author{Tomasz Brzezi\'nski}
\address{Department of Mathematics, University of Wales Swansea,
Singleton Park, Swansea SA2 8PP, U.K.}
\email{T.Brzezinski@swansea.ac.uk}
\urladdr{http//www-maths.swan.ac.uk/staff/tb}
\author{Lars Kadison}
\address{Department of Mathematics and Statistics,  University of New Hampshire, Kingsbury Hall,
Durham, NH 03824} 
\email{kadison@math.unh.edu}
\urladdr{www.math.ntnu.no/\~{}kadison}
\author{Robert Wisbauer}
\address{Department of Mathematics, Heinrich-Heine-University, 
 D-40225
D\"usseldorf, Germany}
\email{wisbauer@math.uni-duesseldorf.de}
\date{June 2002}
\begin{abstract}
A relationship between coseparable corings and separable non-unital 
rings is established. In particular it is shown that an $A$-coring 
$\cC$ has an associative $A$-balanced product. A Morita context is 
constructed for a coseparable coring with a grouplike element. 
Biseparable corings are defined, and a conjecture relating them to 
Frobenius corings is proposed. 
\end{abstract}
\maketitle

\section{Introduction}
Corings were introduced by Sweedler in \cite{Swe:pre} as a 
generalisation of coalgebras and a means for dualising the 
Jacobson-Bourbaki theorem. Recently, corings have resurfaced in the 
theory of Hopf-type modules, in particular it has been shown in 
\cite{Brz:str} that the category of entwined modules 
is an example of a category of comodules of a  coring. Since 
entwined modules appear to be the most general of Hopf-type modules 
studied since the mid-seventies, the theory of corings provides one with 
a uniform and general approach to studying all such modules. This simple 
observation renewed interest in general theory of corings.

Corings appear naturally in the theory of ring extensions. Indeed, they 
provide an equivalent description of certain types of extensions  
(cf.\ \cite{Brz:tow}). In this paper we study properties of corings 
associated to extensions. In particular, we study {\em coseparable 
corings} introduced by Guzman \cite{Guz:coi} (and recently studied in 
\cite{GomLou:cos} from a different point of view) and we reveal an
intriguing duality between such corings and a non-unital 
generalisation of separable ring extensions. We also show that to any 
grouplike element in a coseparable coring one can associate a Morita 
context. This leads to a pair of adjoint functors. One of these 
functors turns out to be fully faithful. Furthermore we introduce the 
notion of a {\em biseparable coring} and study its relationship to 
Frobenius corings introduced in \cite{Brz:tow}. This allows us to 
consider  a conjecture from 
\cite{CaeKad:bis}, concerning biseparable and Frobenius extensions in 
a new framework.

Our paper is organised as follows. In the next section, apart from 
recalling some basic facts about corings and comodules, we introduce a 
non-unital generalisation of separable extensions, which we term 
{\em separable $A$-rings}. We show that any coseparable coring is an 
example of such a separable $A$-ring, and conversely, that every 
separable $A$-ring leads to a non-unital coring. We then proceed in 
Section~3 to construct a Morita context associated to a grouplike 
element in a coseparable coring. We consider some examples coming from
 ring extensions and bialgebroids.  
 Finally in Section~4 we introduce the 
notion of {\em biseparable corings}. These are closely related to 
biseparable extensions, and may serve as a means for settling the
question put forward in \cite{CaeKad:bis} of whether biseparable extensions 
are Frobenius.

Throughout the paper, $A$ denotes an associative ring with unit 
$1_{A}$, and we use the standard notation for right (resp.\ left) $A$-modules 
$\M_{A}$ (resp.\ ${}_{A}\M$), bimodules, such as $\rhom A--$ for 
right $A$-module maps, $\lhom A--$ for left $A$-module maps etc. For any $(A,A)$-bimodule $M$ the centraliser of $A$ in $M$ is denoted by $M^A$, i.e., $M^A := \{ m\in M\; |\; \forall a\in A, \; am=ma\}$.

\section{Coseparable $A$-corings and separable $A$-rings}

\subsection{Coseparable corings}
We begin by recalling the definition of a coring from \cite{Swe:pre}. 
An $(A,A)$-bimodule $\cC$ is said to be a {\em non-counital $A$-coring}
 if there exists an $(A,A)$-bimodule map 
 $\Delta_\cC:\cC\to  \cC\otimes_A\cC$ rendering the following diagram 
 commutative
 $$\xymatrix{
     \cC\ar[rr]^{ {\DC}}\ar[d]_{\DC}  & & \cC\ot_A\cC
\ar[d]^{\id_{\cC}\ut\DC} \\  
 \cC\ot_A\cC \ar[rr]^{\DC\ut\id_{\cC}}  & &\cC\ot_A\cC\ot_A \cC\, . } 
 $$
 The map $\DC$ is termed a {\em
coproduct}. Given a non-counital $A$-coring $\cC$ with a coproduct $\DC$, an  
$(A,A)$-bimodule map $\eps_\cC:\cC\to A$ such that 
$$
(\eps_\cC\ut\id_\cC)\circ\Delta_\cC =  
    \id_{A\otimes_{A}\cC}=\id_\cC, \quad 
  (\id_\cC\ut\eps_\cC)\circ \Delta_\cC =\id_{\cC\otimes_{A}A}=\id_\cC.
  $$
is called a {\em counit} of $\cC$. A non-counital $A$-coring with a 
counit is called an {\em $A$-coring}.

If $\cC$ is an (non-counital) $A$-coring, a right $A$-module $M$ is called 
a {\em non-counital right 
$\cC$-comodule} if there exists a right $A$-module map 
$\varrho^{M}:M\to M\otimes_{A}\cC$ rendering the following diagram 
commutative
$$\xymatrix{
     M\ar[rr]^{{\varrho^M}}\ar[d]_{\varrho^M}  & &M\ot_A\cC
\ar[d]^{\varrho^M\ut\id_{\cC}} \\  
 M\ot_A\cC \ar[rr]^{\id_{M}\ut\DC}  & & M\ot_A \cC\ot_A \cC\,. } 
 $$
The  map $\varrho^M$ is called a {\em $\cC$-coaction}. If, in addition, a 
$\cC$-coaction satisfies  the condition 
$$
(\id_{M}\ut\eC)\circ \varrho^{M} = \id_{M\ot_{A}A} = \id_{M},
$$
then $M$ is called a {\em right $\cC$-comodule}. Similarly one defines 
left $\cC$-comodules, and $(\cC,\cC)$-bicomodules. 
Given right $\cC$-comodules  $M$, $N$, a right 
$A$-linear map $f: M\to N$ is called a {\em morphism of right $\cC$-comodules} 
provided the following diagram
$$\xymatrix{
 M \ar[r]^{f} \ar[d]_{\varrho^M} & N \ar[d]^{\varrho^N} \\
 M\ot_A\cC\ar[r]^{f\ut\id_{\cC}}  & N\ot_A\cC  \,} $$ 
 is commutative. The category of right $\cC$-comodules is denoted 
 by $\M^{\cC}$. We use Sweedler notation to denote the action of 
 a coproduct or a coaction on elements,
 $$
 \DC(c) = \sum c\sw 1\ut c\sw 2, \qquad \varrho^{M}(m) = \sum m\sw 
 0\ut m\sw 1.
 $$
 An immediate example of a left and right $\cC$-comodule is provided 
 by $\cC$ itself. In both cases coaction is given by the coproduct 
 $\DC$. Also, for any right (resp.\ left) $A$-module $M$, the tensor 
 product $M\otimes_{A}\cC$ (resp.\ $\cC\ot_{A}M$)  is a right (resp.\ 
 left) $\cC$-comodule with the coaction $\id_{M}\ut\DC$ (resp.\ 
 $\DC\ut\id_{M}$). This defines a functor which is the right adjoint 
 of a forgetful functor from the category of $\cC$-comodules to the 
 category of $A$-modules. Although this functor can be defined for 
 non-counital corings and non-counital comodules, 
 the adjointness holds only for corings with a 
 counit. 
 
 In particular $\cC\otimes_{A}\cC$ is a 
 $(\cC,\cC)$-bicomodule, and $\DC$ is a $(\cC,\cC)$-bicomodule map, 
 and following \cite{Guz:coi} we have 
\begin{definition} \Label{def.cosep}
   A (non-counital) coring $\cC$ is said to be 
 {\em coseparable} if there exists a $(\cC,\cC)$-bicomodule splitting 
 of the coproduct $\DC$. 
 \end{definition}
 Although Definition~\ref{def.cosep} makes sense for non-counital 
 corings, it is much more meaningful in the case of corings with a 
 counit. In this case $(\cC,\cC)$-bicomodule splittings of $\DC$, 
 $\pi:\cC\ot_{A}\cC\to\cC$ are in bijective correspondence with 
 $(A,A)$-bimodule maps $\gamma:\cC\ot_{A}\cC\to A$ such that for all 
 $c,c'\in\cC$,
 $$
 \sum\gamma(c\ut c'\sw 1)c'\sw 2 = \sum c\sw 1\gamma(c\sw 2\ut c'), 
 \qquad \sum \gamma(c\sw 1\ut c\sw 2) = \eC(c).
 $$
 Such a map $\gamma$ is termed a {\em cointegral} in $\cC$, and the 
 first of the above equations is said to express a {\em colinearity} of 
 a cointegral. The 
 correspondence is given by $\gamma = \eC\circ \pi$ and $\pi(c\ut 
 c') = \sum c\sw 1\gamma(c\sw 2\ut c')$.
 Furthermore $\cC$ is a coseparable $A$-coring if and only if the 
 forgetful functor $\M^{\cC}\to \M_{A}$ is separable (cf.\ 
 \cite[Theorem~3.5]{Brz:str}). 
 
 Corings appear naturally in the context of ring extensions.    A 
 ring extension $B \to A$ determines the \textit{canonical} Sweedler $A$-coring $\cC := 
A \otimes_B A$ with coproduct $\DC: \cC \to \cC \ot_A \cC$
given by $\DC(a \ut a') = a \ut 1_A \ut a'$ and counit $\eC: \cC \to A$
given by $\eC(a \ut a') = aa'$ for all $a,a' \in A$. Recall from 
\cite{Pie:alg} that an extension $B\to A$ is said to be {\em split} 
if there exists a $(B,B)$-bimodule map $E:A\to B$ such that 
$E(1_{A}) = 1_{B}$.  The map $E$ is  known as a {\em conditional 
expectation}. The canonical Sweedler coring associated to a 
split ring extension is coseparable.
A cointegral $\gamma$ coincides with the splitting map $E$ via the 
natural isomorphisms 
$$
\hom BBAB \subset \hom BBAA \cong \hom AA{A\ot_{B}A\ot_{A}A\ot_{B}A}A
$$ 
(cf.\ \cite[Corollary~3.7]{Brz:str}).

 \subsection{Separable $A$-rings}
 Corings can be seen as a dualisation of {\em $A$-rings} and
 coseparable corings turn out to be closely related to a
generalisation of separable extensions of rings. In this subsection 
we describe this generalisation.

\begin{definition}\Label{def.A-ring}
    An $(A,A)$-bimodule  $B$ is called an {\em $A$-ring} provided there 
exists an $(A,A)$-bimodule map $\mu : B\otimes_{A}B\to B$ rendering 
commutative the diagram
 $$\xymatrix{
 B\ot_AB\ot_AB\ar[r]^{\quad\id_{B}\ut \mu}\ar[d]_{\mu\ut\id_{B}} 
 & B\ot_AB \ar[d]^\mu \\ 
  B\ot_AB \ar[r]^{\mu}  & B \,.} $$
 \end{definition}
 This means that $\mu$ is associative. Note that an $A$-ring is necessarily a (non-unital) ring in the usual 
 sense. Equivalently, 
 an $A$-ring can be defined as a ring and an $(A,A)$-bimodule $B$ 
 with product that is an $A$-balanced $(A,A)$-bimodule map.

Note further that the notion of an $A$-ring in 
Definition~\ref{def.A-ring} is a non-unital generalisation of ring 
extensions. Indeed, it is only natural to call an $(A,A)$-bilinear map $\iota:A\to B$ 
a {\em unit} (for $(B,\mu)$)
 if it induces a commutative diagram 
$$\xymatrix{
   B\ar[r]^{\id_{B}\ut\iota\quad}\ar[d]_{\iota\ut\id_{B}}\ar[rd]^= 
   & B\ot_AB\ar[d]^\mu   \\    
 B\ot_AB \ar[r]^\mu & B\, .}$$
 If this holds then $\iota(1_A)=1_B$ is a unit of $B$ in the usual sense. 
 One can then easily show that $\iota$ is a ring map, hence a unital 
 $A$-ring is simply a ring extension.
 
 \begin{definition}\Label{def.mod}
     Given an $A$-ring $B$, a right $A$-module $M$ is said to be a 
     {\em right $B$-module} provided there exists a right $A$-module 
     map $\varrho_M: M\ot_A B\to M$ making the following diagram
     $$\xymatrix{
 M\ot_AB\ot_AB\ar[rr]^{{\varrho_M}\ut\id_{B}}\ar[d]_{\id_{M}\ut\mu}  & 
 & M\ot_AB
\ar[d]^{\varrho_M} \\  
 M\ot_AB \ar[rr]^{\varrho_M}  && M \,} $$
 commute. The map $\varrho_{M}$ is called a {\em right $B$-action}. On 
 elements the action is denoted by a dot in a standard way, i.e., 
 $m\cdot b = \varrho_{M}(m\ut b)$. Remember that for all $a\in A$ we 
 have $(ma)\cdot b = m\cdot (ab)$.
 
 A morphism $f:M\to N$ between two $B$-modules is 
 an $A$-linear map which makes the following diagram
  $$\xymatrix{
 M\ot_AB \ar[r]^{f\ut\id_{B}}\ar[d]_{\varrho_M}  & N\ot_AB \ar[d]^{\varrho_N} \\ 
  M \ar[r]^{f}  & N \,} $$
 commute.  A right $B$-module $M$ is said to 
 be {\em firm} provided the induced map $M\otimes_{B}B\to M$, $m\ut 
 b\mapsto m\cdot b$ is a right $B$-module isomorphism. The category 
 of firm right $B$-modules is denoted by $\M_{B}$. Note that, by definition, $\M_B$ is a subcategory of the category of right $A$-modules $\M_A$.
 \end{definition}
 
 Obviously left $B$-modules are defined in a symmetric way. Similarly 
 one defines $(B,B)$-bimodules, $(A,B)$-bimodules etc. 
 
 Dually to the definition of coseparable corings we can define 
 separable $A$-rings.
 \begin{definition}\Label{def.sep.A-ring}
     An  $A$-ring $B$ is said to be {\em separable} 
 if the product map $\mu:B\ot_AB\to B$ has a $(B,B)$-bimodule  section 
 $\delta:B\to B\ot_AB$. 
\end{definition}
If $B$ is a separable $A$-ring then clearly $\mu$ is surjective and the 
induced map  $B\otimes_{B}B\to B$ is an isomorphism. Therefore $B$ is 
a firm left and right $B$-module,  i.e., $B$ is a firm ring. 

Note that if $B$ has a unit $\iota:A\to B$ then $B$ is a
separable $A$-ring if and only if $B$ is a separable extension of 
$A$. Thus Definition~\ref{def.sep.A-ring} extends the notion of a separable 
extension to non-unital rings. Note, however, that in general this is not an extension, 
since there is no (ring) map $A\to B$.

\begin{remark}
    In consistency with $A$-corings,   we use the terminology of 
 \cite{BerHau:cog} in Definition~\ref{def.A-ring}.  In \cite[11.7]{Pie:alg} $A$-rings are termed
 {\em multiplicative $A$-bimodules}. 
Following \cite{Tay:big} one might call a
separable $A$-ring (as defined in Definition~\ref{def.sep.A-ring}) an {\em $A$-ring with a
splitting map}.
\end{remark}

\subsection{Coseparable $A$-corings are separable $A$-rings}
The main result of this section is contained in the following

\begin{theorem}
    If $\cC$ is a coseparable $A$-coring then  $\cC$ is a separable 
    $A$-ring.
\label{sep.cosep.thm}
\end{theorem}
\begin{proof}
  Let $\pi: \cC\otimes_{A}\cC\to \cC$ be a bicomodule retraction of 
    the coproduct $\DC$, and let $\gamma = \eC\circ \pi$ be 
    the corresponding cointegral. We claim that  $\cC$ is an associative 
    $A$-ring with product $\mu = \pi$. Indeed,
    since
the alternative expressions for product 
are $cc' =  \sum \gamma(c\ut c'\sw 1)c'\sw 2 = \sum c\sw 
1\gamma(c\sw 2\ut c')$,  
 for all $c,c',c''\in\cC$ we have, 
using the left $A$-linearity of $\gamma$ and $\DC$, 
$$
(cc')c'' = \sum (\gamma(c\ut c'\sw 1)c'\sw 2)c'' = 
 \sum \gamma(c\ut c'\sw 1)\gamma(c'\sw 2 \ut c''\sw 1)c''\sw 2.
$$
On the other hand, the colinearity and 
right $A$-linearity of $\gamma$, and the left $A$-linearity of $\DC$ imply
\begin{eqnarray*}
    c(c'c'') & = & \sum c(\gamma(c'\ut c''\sw 1)c''\sw 2) = 
    \sum \gamma(c\ut \gamma(c'\ut c''\sw 1)c''\sw 2)c''\sw 3 
    \\
    &=& \sum \gamma(c\ut c'\sw 1\gamma(c'\sw 2 \ut c''\sw 1))c''\sw 2
    = \sum \gamma(c\ut c'\sw 1)\gamma(c'\sw 2 \ut c''\sw 
    1)c''\sw 2. 
\end{eqnarray*}
This explicitly proves that the product in $\cC$ is associative. 
Clearly this product is $(A,A)$-bilinear. Note that $\DC$ is a $(\cC,\cC)$-bimodule 
map since 
$$
c\DC(c') = \sum cc'\sw 1\ut c'\sw 2 = \sum \pi(c\ut c'\sw 1)\ut 
c'\sw 2 = \DC\circ \pi(c\ut c') = \DC(cc'),
$$
from the right colinearity of $\pi$. Similarly for the left 
$\cC$-linearity. Finally $\pi$ is split by $\DC$ since $\pi$ is a 
retraction of $\DC$. This proves that $\cC$ is a separable $A$-ring. 
\end{proof}

\begin{proposition}
    Let $\cC$ be a coseparable $A$-coring with a cointegral $\gamma$. 
    View $\cC$ as an $A$-ring with product $\pi$ as in 
    Theorem~\ref{sep.cosep.thm}. 
    Then any right $\cC$-comodule $M$ 
    is a firm right $\cC$-module with the action $\varrho_{M} = 
    (\id_{M}\ut \gamma)\circ (\rM\ut \id_{\cC})$. 
\label{cosep.firm.mod}
\end{proposition}
\begin{proof}
    Take any $m\in M$ and $c,c'\in \cC$. Then explicitly the action 
    reads  $m\cdot c = \sum m\sw 0\gamma(m\sw 
    1\ut c)$, and  we can compute
    \begin{eqnarray*}
	(m\cdot c)\cdot c' &=& \sum (m\sw 0\gamma(m\sw 1\ut c))\cdot c' =\sum m\sw 0\gamma(m\sw 
	1\gamma(m\sw 2\ut c)\ut c') \\
	&=&\sum m\sw 0\gamma(\gamma(m\sw 1\ut c\sw 1)c\sw 2\ut c')= 
	\sum m\sw 0\gamma(m\sw 1\ut c\sw 1)\gamma(c\sw 2\ut c')\\
	&=& \sum m\sw 0\gamma(m\sw 1\ut c\sw 1\gamma(c\sw 2\ut c')) = \sum m\sw 
	0\gamma(m\sw 1\ut cc') = m\cdot (cc'),
     \end{eqnarray*}
     as required. We used the following properties of a cointegral: 
     colinearity to derive the third equality and $A$-bilinearity to 
     derive the fourth and fifth equalities. Obviously the 
     action is right $A$-linear. Thus $M$ is a $\cC$-module. We need 
     to show that it is firm.
     
     Note that $M\otimes_{\cC}\cC$ is defined as a cokernel of the  
     following right $\cC$-linear map
     $$
     \lambda: M\otimes_{A}\cC\otimes_{A}\cC \to M\otimes_{A}\cC, 
     \quad m\ut c\ut c' \mapsto mc\ut c' - m\ut cc'.
     $$
     Since $\gamma$ is a cointegral, $\varrho_{M}$ is a right 
     $\cC$-linear retraction of $\varrho^{M}$, hence, in particular it is 
     a surjection and we have the following sequence of right 
     $\cC$-module maps
     $$
     \xymatrix@1{
     M\otimes_{A}\cC\otimes_{A}\cC \ar[r]^-{\lambda} & M\otimes_{A}\cC  
     \ar[r]^-{\varrho_{M}}&
     M \ar[r] & 0.}
     $$
     We need to show that this sequence is exact. Clearly the 
     associativity of action of $\cC$ on  $M$ implies that 
     $\varrho_{M}\circ\lambda =0$, so that $\im\lambda\subseteq 
     \ker\varrho_{M}$. Furthermore, for all $m\in M$ and 
     $c\in \cC$ we have
     \begin{eqnarray*}
	 (\rM\circ\varrho_{M}-\lambda\circ(\id_{M}\ut\DC))(m\ut c)
	 &=& \sum m\sw 0\ut m\sw 1\gamma(m\sw 2\ut c)\\
	 &&- \sum m\cdot c\sw 1\ut
	 c\sw 2
	 + \sum m\ut \pi(c\sw 1\ut c\sw 2)\\
	 &=& \sum m\sw 0 \gamma(m\sw 1\ut c\sw 1)\ut c\sw 2\\
	 &&- 
	 \sum m\sw 0 \gamma(m\sw 1\ut c\sw 1)\ut c\sw 2 + 
	 m\ut c\\
	 &=&m\ut c,
     \end{eqnarray*}
     where we used the colinearity of a cointegral. This implies that 
     $\ker\varrho_{M}\subseteq 
     \im\lambda$, i.e., the above sequence is exact as required.
\end{proof}

As an example of a coseparable coring one can take the canonical 
Sweedler coring associated to a split ring extension $B\to A$. In 
this case the product in $A\ot_{B}A$ comes out as
$$
(a\ut a')(a''\ut a''') = aE(a'a'')\ut a''', \qquad \forall 
a,a',a'',a'''\in A
$$
where $E$ is a splitting map. This is known as the {\em 
$E$-multiplication}. Since a comodule of the canonical coring is a 
{\em descent datum} for a ring extension $B\to A$, 
Proposition~\ref{cosep.firm.mod} implies that every descent datum is 
a firm module of the $A$-ring $A\ot_{B}A$ with the $E$-multiplication.

Theorem~\ref{sep.cosep.thm} has the following (part-) converse.

\begin{proposition}\Label{sep.A-ring.cor}
Let $B$ be a separable $A$-ring. Then $B$ is a coseparable 
non-counital coring.
\end{proposition}
\begin{proof}   
  Let   $\Delta: B\to B\otimes_A B$ be 
 a $(B,B)$-bimodule map  splitting 
the product $\mu$ in $B$.  The $B$-linearity of  $\Delta$ implies that 
the following diagram
$$\xymatrix{  
B\ot_A B\ot_A B \ar[d]_{ \mu\ot\id_{B}}&& \ar[ll]_{\quad\id_{B}\ot\Delta} B\ot_A B
 \ar[d]^{\mu}  \ar[rr]^{\Delta\ot \id_{B}\quad}  && B\ot_A B\ot_A B 
 \ar[d]^{\id_{B}\ot\mu} \\
 B\ot_A B &&\ar[ll]^{\quad\Delta}  B \ar[rr]_{\Delta\quad} && B\ot_A B \,}$$
is commutative. For all $b\in B$ we write $(\Delta\ut \id_{B})\circ \Delta(b) 
= \sum b\sw 1\sw 
1\ut b\sw 1\sw 2\ut b\sw 2$ and  $(\id_{B}\ut\Delta)\circ \Delta(b) 
      = \sum b\sw{1} \ut b\sw{2}\sw{1}\ut b\sw{2}\sw{2}$, and use the 
      above diagram to obtain 
$$\begin{array}{rl}
 \Delta(b)=(\Delta\circ\mu\circ \Delta)(b)& = 
 (\id_{B}\ut\mu)\circ (\Delta\ut \id_{B})\circ \Delta(b)
       = \sum \; b\sw{1}\sw{1} \ut \mu (b\sw{1}\sw{2}\ut b\sw{2}) \\
   &=  (\mu\ut\id_{B})\circ (\id_{B}\ut\Delta)\circ \Delta(b) 
      = \sum \;\mu(b\sw{1} \ut b\sw{2}\sw{1})\ut b\sw{2}\sw{2}.
 \end{array}$$
Using these identities  we can compute
$$\begin{array}{rl}
 (\id_{B}\ut\Delta)\circ \Delta(b) &= 
\sum b\sw{1}\sw{1}\ut (\Delta\circ \mu)(b\sw{1}\sw{2}\ut b\sw{2}) \\
 &= \sum b\sw{1}\sw{1}\ut  ((\mu\ut\id_{B})\circ(\id_{B}\ut\Delta)) 
 (b\sw{1}\sw{2}\ut b\sw{2}) \\
&= \sum b\sw{1}\sw{1}\ut \mu(b\sw{1}\sw{2}\ut b\sw{2}\sw{1})\ut 
b\sw{2}\sw{2}\, ,
\end{array}
$$
and 
  $$
  \begin{array}{rl}
 (\Delta\ut \id_{B})\circ \Delta(b) &= 
  \sum (\Delta\circ\mu) (b\sw{1}\ut b\sw{2}\sw{1})\ut b\sw{2}\sw{2} \\
&= \sum ((\id_{B}\ut\mu)\circ(\Delta\ut\id_{B}))(b\sw{1}\ut 
b\sw{2}\sw{1})\ut 
b\sw{2}\sw{2}\\
&=\sum b\sw{1}\sw{1}\ut \mu (b\sw{1}\sw{2}\ut b\sw{2}\sw{1})\ut 
b\sw{2}\sw{2}\, ,
\end{array} $$
i.e., $(\Delta\ut \id_{B})\circ \Delta = (\id_{B}\ut\Delta)\circ \Delta$. 
This proves that $B$ is a non-counital $A$-coring with coproduct 
$\Delta$.

Next note 
that the above diagram can also be 
understood as a statement that $\mu$ is a $(B,B)$-bicomodule 
map. Since $\mu$ is a retraction for $\Delta$,  $B$ is a coseparable 
non-counital coring as required. This completes the proof.
\end{proof}

\section{Morita contexts for coseparable corings}
Although the Morita theory is usually developed for rings with unit, 
it can be extended to firm rings without units (cf.\ 
\cite[Exercise~4.1.4]{BerKea:cat}). Recall that a right module $M$ of 
a non-unital ring $R$ 
is said to be {\em firm} 
if the map $M\otimes_{R}R\to M$ induced from the $R$-product in $M$ is 
an $R$-module isomorphism. Similarly one defines left firm modules. A 
non-unital
ring is a firm ring if it is firm as a left and right $R$-module.

\begin{definition}\label{def.Mor.firm}
    Given a pair  of firm non-unital rings $R,S$, a Morita context 
    consists of a firm $(R,S)$-bimodule $V$ and a firm 
    $(S,R)$-bimodule $W$ and a pair of bimodule maps
    $$
    \sigma : W\otimes_{R}V\to S, \qquad \tau :
    V\otimes_{S}W\to R,
    $$
    such that the following diagrams
    $$
    \xymatrix{ W\ot_{R}V\ot_{S}W 
    \ar[rr]^{\sigma\ut\id_{W}}\ar[d]_{\id_{W}\ut\tau} && S\ot_{S}W 
    \ar[d]^{\cong}& & V\ot_{S}W\ot_{R}V
    \ar[rr]^{\tau\ut\id_{V}}\ar[d]_{\id_{V}\ut\sigma} && R\ot_{R}V
    \ar[d]^{\cong}\\
    W\ot_{R}R \ar[rr]^{\cong} && W && V\ot_{S}S \ar[rr]^{\cong} && V}
    $$
    commute. A Morita context is denoted by $(R,S,V,W,\tau,\sigma)$. A Morita 
    context is said to be {\em strict} provided $\sigma$ and $\tau$ 
    are isomorphisms.
\end{definition}

The Morita theory for non-unital rings can be developed along the same 
lines as the usual Morita theory \cite{Q:notes}. The aim of this section is to show 
that there is a Morita context associated to any coseparable coring 
with a grouplike element. To construct such a context we employ 
techniques developed in recent papers \cite{Abu:Mor}, \cite{Cae:Mor}.

First recall that an element $g$ of an $A$-coring $\cC$ is said to be 
a {\em grouplike element}, provided $\DC(g) = g\ut g$ and $\eC(g) =1$. 
Obviously, not every coring has grouplike elements. The results of 
this section are contained in the following
\begin{theorem} \label{thm.Morita}
    Let $\cC$ be a coseparable $A$-coring with a cointegral $\gamma$ 
    and a grouplike element $g$. 
    View $\cC$ as a separable $A$-ring as in 
    Theorem~\ref{sep.cosep.thm} and let $\M_{\cC}$ denote
    the category of firm right modules of the $A$-ring $\cC$ (cf.\ 
    Definition~\ref{def.mod}). For any $M\in \M_{\cC}$ define
    $$
    M_{g,\gamma}^{\cC} = \{m\in M \; |\; \forall c\in \cC, \; m\cdot 
    c = m\gamma(g\ut c)\}.
    $$
    Then:
    \begin{zlist}
	\item $B= A_{g,\gamma}^{\cC} = \{b\in A \; |\; \forall c\in \cC, \; 
	\gamma(gb\ut c) = b\gamma(g\ut c)\}$,
	is a subring of $A$.
	\item The assignment $( -)_{g,\gamma}^{\cC}: \M_{\cC}\to \M_{B}$, 
	$M\mapsto M_{g,\gamma}^{\cC}$ is a covariant functor which has a left 
	adjoint $-\otimes_{B}A:\M_{B}\to\M_{\cC}$. 
	\item $Q=\cC_{g,\gamma}^{\cC}$ is a firm left ideal in $\cC$ and 
	hence a $(\cC,B)$-bimodule.
	\item For every $M\in \M_{\cC}$, the additive map
	$$
	\omega_{M}:M\otimes_{\cC}Q\to M_{g,\gamma}^{\cC}, \qquad m\ut q\mapsto m\cdot 
	q,
	$$
	is bijective.
	\item Define  two maps 
	$$
	\sigma :Q\ot_{B}A\to \cC, \;q\ut a\mapsto qa \;\; 
	\mbox{\rm and} \;\;
	\tau: A\ot_{\cC}Q\to B, \; a\ut q\mapsto \gamma(ga\ut q). 
	$$
	Then 
	$(B,\cC,A,Q,\tau,\sigma)$ is a Morita context, in which $\tau$ is 
	surjective (hence an isomorphism).
	\end{zlist}
\end{theorem}
\begin{proof}
    
    (1) First note that $A$ is a right $\cC$-comodule with the 
    coaction $a\mapsto 1_{A}\ut ga$. Thus by 
    Proposition~\ref{cosep.firm.mod}, $A$ is a firm $\cC$-module with 
    the action $a\cdot c = \gamma(ga\ut c)$, for all $a\in A$ and 
    $c\in \cC$. Therefore the definition of $B$ makes sense and takes 
    the form stated. Obviously $1_{A}\in B$. Furthermore for all $b,b'\in 
    B$ and $c\in \cC$ we have
    $$
    \gamma(gbb'\ut c) = \gamma(gb\ut b'c) = b\gamma(g\ut b'c) = 
    bb'\gamma(g\ut c),
    $$ 
    so that $bb'\in B$ as required. Alternatively, we note that $B \cong \rend \cC A$
and is therefore a ring. 
    
    (2) We note that $M^{\cC}_{g,\gamma} \cong \rhom \cC AM$ via $f \mapsto f(1_A)$
has the left adjoint $-\ot_B A$.  In more detail, take any $M\in \M_{\cC}$, 
$m\in M_{g,\gamma}^{\cC}$, $b\in B$ and $c\in \cC$, and 
    use the definitions of $B$ and $M_{g,\gamma}^{\cC}$ to compute
    $$
    (mb)\cdot c = m\cdot (bc) = m\gamma(g\ut bc) = mb\gamma(g\ut c).
    $$
    This shows that $M_{g,\gamma}^{\cC}$ is a right $B$-module,
    hence $(-)_{g,\gamma}^{\cC}$ is a functor as stated. In the 
    opposite direction, for any right $B$-module $N$, $N\ot_{B}A$ 
    is a firm right $\cC$-module 
    with the action $(n\ut a)\cdot c = n \ut \gamma(ga\ut c)$. Note 
    that this action is well-defined by the construction of $B$. More 
    precisely,  $A$ is a left $B$-module and a firm right 
     $\cC$-module (since it is a right $\cC$-comodule). It is a 
     $(B,\cC)$-bimodule, since for every $b\in B$, $a\in A$ and $c\in 
     \cC$ we have
     $$
     (ba)\cdot c = \gamma(gba\ut c) = \gamma (gb\ut ac) = 
     b\gamma(g\ut ac) = b(a\cdot c),
     $$
     by (1). The above action is simply induced from the action of 
     $\cC$ on $A$ and thus well-defined. Therefore
    we have a functor as required. Now, one can easily check that the 
    unit and counit of the adjunction are given by
    $$
    \eta_{N}: N\to (N\ut_{B}A)_{g,\gamma}^{\cC}, \qquad n\mapsto n\ut 
    1_{A},
    $$
    $$
    \eps_{M}: M_{g,\gamma}^{\cC}\ot_{B}A\to M, \qquad m\ut a\mapsto 
    ma,
    $$ 
    for all $N\in \M_{B}$ and $M\in \M_{\cC}$.
    
    (3) Note that $Q \cong \rhom \cC A \cC$ and is therefore a natural $(\cC,B)$-bimodule.  In more detail, we have
    $(cq)c' = c(qc') = cq\gamma(g\ut c')$
 for any $c, c'\in \cC$ and $q\in Q$ , so $cq$ is an element of $Q$, 
    hence $Q$ is a left $\cC$-ideal. By (2) $Q$ is also a right 
    $B$-module, and since the product in $\cC$ is an $(A,A)$-bimodule 
    map, $Q$ is a $(\cC,B)$-bimodule. We only need to show that $Q$ 
    is firm as a left $\cC$-module. This can be shown by the same 
    technique as in the proof of Proposition~\ref{cosep.firm.mod}. 
    $\cC\otimes_{\cC}Q$ is defined as a cokernel of the  
     following left $\cC$-linear map
     $$
     \lambda: \cC\otimes_{A}\cC\otimes_{A}Q \to \cC\otimes_{A}Q, 
     \quad \lambda = \id_{\cC}\ut \pi - \pi\ut \id_{Q},
     $$
     where $\pi$ is the product map in $\cC$ (i.e., the splitting of 
     $\DC$) corresponding to the 
     cointegral $\gamma$. Observe that the product map $ \pi : \cC\ot_{A}Q\to Q$ is a 
     surjection. Indeed, first note that since for all $c\in\cC$, 
     $\pi(g\ut c) = g\gamma(g\ut c)$ by the relationship between 
     $\pi$ and $\gamma$, the grouplike element $g$ is in $Q$. For 
     any $q\in Q$ take $q\ut g\in \cC\ot_{A}Q$. Then $\pi(q\ut g) = 
     q\gamma(g\ut g) = q$, by the properties of the cointegral $\gamma$. 
     Thus $\pi$ is a surjection as claimed. 
     
     We have the following sequence of left
     $\cC$-module maps
     $$
     \xymatrix@1{
     \cC\otimes_{A}\cC\otimes_{A}Q \ar[r]^-{\lambda} & \cC\otimes_{A} Q 
     \ar[r]^-{\pi}&
     Q \ar[r] & 0.}
     $$
     We need to show that this sequence is exact. Clearly the 
     associativity  of $\pi$ implies that 
     $\pi\circ\lambda =0$, so that $\im\lambda\subseteq 
     \ker\pi$. Furthermore,  we have
     $$
     \DC\circ\pi-\lambda\circ(\DC\ut\id_{Q}) = 
     (\id_{\cC}\ut\pi)\circ (\DC\ut\id_{Q}) -(\id_{\cC}\ut\pi)\circ 
     (\DC\ut\id_{Q}) + (\pi\circ\DC)\ut\id_{Q} = \id_{\cC}\ut\id_{Q},
     $$
     where we used the colinearity of $\pi$ and the fact that $\pi$ is 
     a splitting of $\DC$. This implies that 
     $\ker\pi\subseteq 
     \im\lambda$, i.e., the above sequence is exact as required.
     
     (4) Note that $\omega_M$ is the natural map $M \ot_{\cC} \rhom \cC A \cC
\to \rhom \cC A M \cong M^{\cC}_{g,\gamma}$, $m\ot f\mapsto [a\mapsto f(1)m]$, and hence it is well-defined.  Explicitly, 
 for all 
     $m\in M$, $q\in Q$ and $c\in\cC$ we have $(m\cdot q)\cdot c = 
     m\cdot (qc) = m\cdot q\gamma(g\ut c)$, as required.  We need to show that $\omega_M$ is bijective. Consider a map
     $$
     \theta_{M}:M_{g,\gamma}^{\cC} \to M\ot_{\cC}Q, \qquad m\mapsto 
     m\ut g.
     $$
     This is well-defined, since as shown in the proof of (3), $g\in 
     Q$. Take any $m\in M_{g,\gamma}^{\cC}$. Then
     $$
     \omega_{M}(\theta_{M}(m)) = m\cdot g = m\gamma(g\ut g) = m,
     $$
     by the definition of $M_{g,\gamma}^{\cC}$ and properties of a cointegral. 
     Conversely, for any simple tensor $m\ut q\in M\ot_{\cC}Q$ we have
     $$
     \theta_{M}(\omega_{M}(m\ut q)) = m\cdot q\ut g = m\ut qg = m\ut 
     q\gamma(g\ut g) = m\ut q,
     $$
     again by the definition of $Q$ and properties of cointegrals.
     
     (5) Note that $\sigma$ is the evaluation mapping $\rhom \cC A \cC \ot_B A
\rightarrow \cC$, while $\tau$ is the canonical map $A \ot_{\cC} \rhom \cC A \cC 
\rightarrow \rend \cC A$. 
In more detail, we show that the maps $\sigma$ and $\tau$ are 
     well-defined as bimodule maps. Obviously, $\sigma$ is left 
     $\cC$-linear. Take any $q\in Q$, $a\in 
     A$ and $c\in \cC$ and compute
     $$
     \sigma(q\ut a\cdot c) = q\gamma(ga\ut c) = q\gamma(g\ut ac) = 
     \pi(q\ut ac) = \pi(qa\ut c) = (qa)c.
     $$
     This shows that $\sigma$ is $(\cC,\cC)$-bilinear as required.
     Note that $\tau = \omega_{A}$, and since $B = A_{g,\gamma}^{\cC}$ 
     it is well defined and surjective. 
    Clearly $\tau$ is 
    right $B$-linear. An easy computation which involves the 
    definition of $B$ confirms that $\tau$ is $(B,B)$-bilinear. Next  
    we need to check the commutativity of diagrams in 
    Definition~\ref{def.Mor.firm}. The commutativity of the second 
    diagram follows immediately from $A$-linearity of the cointegral. 
    Now take any $a\in A$, and $q,q'\in Q$ and compute
    $$
    \sigma(q\ut a)q' = (qa)q' = q(aq') = q\gamma(g\ut aq') = 
    q\tau(a\ut q'),
    $$
    where the definition of $Q$ was used to derive the third equality. Thus 
    we have a Morita context as required. A standard argument in 
    Morita theory confirms that $\tau$ is an isomorphism (cf.\ 
    \cite[II~(3.4)~Theorem]{Bas:KTh}).
\end{proof}
\begin{corollary}
    With the assumptions and notation as in Theorem~\ref{thm.Morita} 
    we have:
    \begin{zlist}
	\item $Q$ is a subring of $\cC$ with a right unit $g$.
	\item The functor $-\otimes_{B}A: \M_{B}\to \M_{\cC}$ is fully 
	faithful, i.e., the unit of adjunction $\eta$ is an isomorphism.
	\item $A_\cC$ and ${}_\cC Q$ are direct summands of $\cC^n$ for some (finite) $n\in\mathbb{N}$.
	\item $A$ and $Q$ are generators as left resp.\ right $B$-modules.
    \end{zlist}
 \end{corollary}
\begin{proof} (1) follows immediately from the proof of 
Theorem~\ref{thm.Morita}(3), while 
 (2)-(4) follow from Morita theory with surjective $\tau$, and can be proven by the same methods as in the unital case (cf.\ \cite[Ch.~II.3]{Bas:KTh}). In particular, in the case of (3) the ``dual bases" of $A$ and $Q$ can be constructed as follows. Let $\{ a_i\in A, q^i\in Q\}_{i=1,\ldots , n}$ be such that $1_B = \sum_i\tau(a_i\ut q^i)$. Define $\sigma^i = \sigma(q^i\ut -)\in \rhom \cC A\cC$. Then for every $a\in A$ we have
$$
\sum_i a_i\cdot\sigma^i(a) = \sum_i a_i\cdot \sigma(q^i\ut a) = \sum_i \tau(a_i\ut q^i)a =a,
$$
so that $\{ a_i, \sigma^i\}_{i=1,\ldots , n}$ is a dual basis for $A_\cC$. Similarly a dual basis for ${}_\cC Q$ can be constructed as $\{ q^i, \sigma_i\}$ with $\sigma_i = \sigma(-\ut a_i)$.
 \end{proof}

 One can easily find a sufficient condition for the Morita context of 
 Theorem~\ref{thm.Morita} to be strict.
 \begin{proposition}\label{prop.strict}
     Let $\cC$ be a coseparable $A$-coring with a grouplike 
    element $g$. Then the Morita context $(B,\cC,A,Q,\tau,\sigma)$ 
    associated in Theorem~\ref{thm.Morita} is strict, provided $g$ is 
    a left unit in $\cC$.
\end{proposition}
\begin{proof}
    In this case $\gamma(g\ut c) = \eC(c)$, hence $B$ and $Q$ are 
    characterised by relations $\eC(bc) = b\eC(c)$ and $qc = q\eC(c)$, 
    respectively, for all $c\in \cC$. We need to show that $\sigma$ is 
    an isomorphism. For any $c\in \cC$, $\sigma (g\ut \eC(c)) = 
    g\eC(c) = gc = c$, since $g\in Q$. Thus $\sigma$ is surjective. 
    Suppose now that $\sum_{i}q_{i}\ut a_{i}\in \ker\sigma$, i.e., 
    $\sum_{i}q_{i}a_{i}=0$. This implies that 
    $\sum_{i}\eC(q_{i})a_{i}=0$, so that 
    $$
    0 = g\ut_{B}\sum_{i}\eC(q_{i})a_{i} = \sum_{i}g\eC(q_{i})\ut 
    a_{i} = \sum_{i}gq_{i}\ut 
    a_{i} = \sum_{i}q_{i}\ut 
    a_{i} .
    $$ Here we used that for all $q\in Q$, $\eC(q) = \gamma(g\ut 
    q)\in B$, the fact that $g\in Q$ and that $g$ is a left unit in 
    $\cC$. This completes the proof.
\end{proof}

Finally, we consider two examples of Theorem~\ref{thm.Morita}.
\begin{example}
\label{ex:split}
Consider a split extension $\bar{B} \stackrel{\iota}{\longrightarrow} A$ with 
 splitting map $E :A\to \bar{B}$. Then the canonical Sweedler 
coring $\cC = A\ot_{\bar{B}}A$ is an $A$-ring with the 
$E$-multiplication,  $1_{A}\ut 1_{A}\in \cC$ is a grouplike element, 
and the Morita context constructed in 
Theorem~\ref{thm.Morita} comes out as follows. The ring $B$ is just
$$
B = \iota(\bar{B}),
$$
while the $(\cC,B)$-bimodule $Q \subset A \ot_B A$ is
$$
Q \cong A
$$
via $a \mapsto a \ot 1_A$ and the left module action of the $A$-ring $\cC$ on $A$  given
by $c \cdot a = \sum c^1 E(c^2 a)$ (suppressing a possible summation in
$c =\sum  c^1\ut c^2\in A \ot_B A$). The module $A_{\cC}$ is similarly given by
$a \cdot c = \sum E(ac^1)c^2$ for each $a \in A, c \in \cC$. 
 The Morita maps read $\sigma(a\ut a') = a \ut a'$ and $\tau(a\ut a') = 
 E(aa')$. This context is obviously strict.

The proof of this involves applying the theorem, noting that
$$ B =  
\{ b\in A\; |\; \forall a\in A, \; 1_A  E(ba) = bE(a)\} = \iota(\bar{B})\} $$
since $\supseteq$ is clear and $\subseteq$ follows from letting $a = 1_A$.
Next one notes that
$$ Q = \{ 
q=\sum_{i}q_{i}\ut \bar{q}_{i}\in A\ot_{\bar{B}}A\; |\; \forall a\in A, 
\; \sum_{i} q_{i}E(\bar{q}_{i}a)\ut 1_A = qE(a)\} = A \ut 1_A $$
since $\supseteq$ is clear and $\subseteq$ follows from taking $a = 1_A$.      
\end{example}

\begin{example}  
Hopf algebroids over a noncommutative base $k$-algebra $A$, where $k$ is a commutative ring, provide examples of $A$-corings
with grouplike elements; in particular, the canonical bialgebroids
$\End_k A$ and $A \ot_k A^{\rm op}$ do (cf.\  \cite{Lu:Hop} \cite{BrzMil:bia} for the definition and examples of Hopf algebroids).  They can be extended to ring extensions via an
algebraic formulation of depth two for subfactors \cite{KadSzlach:depth}:
a ring extension $B \rightarrow A$ is of \textit{depth two} (D2) if 
$A \ot_B A$  is isomorphic 
 to a direct summand of $A \oplus \cdots \oplus A$ as a $(B,A)$-bimodule, and similarly
as an $(A,B)$-bimodule.  The two conditions are equivalent respectively 
to  the existence of finitely many elements $c_j,b_i \in (A \ot_B A)^B$ and $\gamma_j, \beta_i \in \lrend BBA$
such that  for all $a,a' \in A$,
\begin{equation}
\label{eqs:depthtwo}
 a \ut a' = \sum_i b_i \beta_i(a)a' = \sum_j a \gamma_j(a')c_j.
\end{equation}

Denoting the centraliser $A^B$ of a D2 extension $B \rightarrow A$    
by $R$, the following $R$-coring structure for $\cC := \lrend BBA$ is 
considered in \cite{KadSzlach:depth}.  The $(R,R)$-bimodule structure is
$r  \alpha  r' = r \alpha(-) r'$ ($\alpha \in \cC$).  The coproduct is given most simply
by noting \cite[3.10]{KadSzlach:depth}: $\cC \ot_R \cC \cong \hom BB{A \ot_B A}A$ via $$\alpha \ut \beta \longmapsto
[a \ut a' \mapsto \alpha(a)\beta(a')].$$  Then 
$$\DC(\alpha)(a \ut a') = \alpha(aa'),$$
with counit
$$ \eC(\alpha) = \alpha(1_A).$$
We also have the alternative formulae for the coproduct \cite[Eqs.\ (66), (68)]{KadSzlach:depth}:
$$ \DC(\alpha) = \sum_j \gamma_j \ot_R c^1_j \alpha(c^2_j -) = \sum_i \alpha(- b^1_i)b^2_i \ot_R \beta_i, $$
where $c_j = \sum c_j^1\ut c_j^2$ and $b_i = \sum b_i^1\ut b_i^2$ is a notation suppressing a possible
 summation index.
We note the grouplike element $\id_A$. 

Suppose the D2 extension $B \rightarrow A$ is separable with separability element
$e = \sum e^1 \ut e^2 \in A \ot_B A$ (summation index suppressed).  
Then the $R$-coring $\cC$ is coseparable with cointegral $\gamma: \cC \ot_R \cC \rightarrow
R$ given by  $\gamma(\alpha \ut \beta) =\sum \alpha(e^1)\beta(e^2)$.  
The corresponding $R$-ring structure on $\cC$ is given by ($x \in A$) 
$$ (\alpha * \beta)(x) = \sum \alpha(xe^1)\beta(e^2) =\sum \alpha(e^1)\beta(e^2 x) $$
with the $R$-bimodule structure above.

The Morita context in the theorem  applied to $\cC$  turns out as follows:
\end{example}

\begin{proposition}
The centre $Z$ of $A$ and the non-unital ring $(\cC,*)$ are related by the Morita context
$(Z,\cC, {}_ZR_{\cC}, {}_{\cC}R_Z,\tau,\sigma)$ where
${}_ZR_{\cC}$ is given by $z  r \cdot \alpha = \sum ze^1 r \alpha(e^2)$,
${}_{\cC}R_Z$ by $\alpha \cdot r z = \sum \alpha(e^1)r e^2 z$,
$$ \sigma: R \ot_Z R \rightarrow \cC,\ \ r \ut r' \longmapsto \lambda(r) \rho(r')
= \rho(r')\lambda(r), $$
where $\lambda(r),\rho(r) \in \cC$ denote left and right multiplication by $r \in R$,
respectively, and
$$\tau: R \ot_{\cC} R \stackrel{\cong}{\longrightarrow} \cC,\ \ r \ut r' \longmapsto
\sum e^1 rr' e^2. $$
The Morita context is strict if $B \rightarrow A$ is H-separable.
\end{proposition}
\begin{proof}
We check that $\gamma$ is a cointegral. Take any $\alpha\in\cC$ and compute
$$ \sum \gamma(\alpha\sw 1 \ut \alpha\sw 2) = \sum \alpha(e^1e^2) = \eC(\alpha). $$
Thus $\gamma$ is normalised. Furthermore, for all $\alpha,\beta\in\cC$,
$$ \sum \gamma(\alpha \ut \beta\sw 1) \beta\sw 2 = \sum_i 
\alpha(e^1) \beta(e^2 b^1_i)b^2_i \beta_i = \sum \alpha(e^1)\beta(e^2-). $$
On the other hand
$$\sum  \alpha\sw 1 \gamma(\alpha\sw 2 \ut \beta) = \sum_j \gamma_j(-) c^1_j \alpha(c^2_j e^1) 
\beta(e^2) =\sum  \alpha(-e^1)\beta(e^2), $$
so that $\gamma$ is colinear and hence a cointegral.

The subring of $R$ in Theorem~\ref{thm.Morita} is
$$R^{\cC}_{\id_A, \gamma} = \{ r \in R \, | \, \forall a \in A, \, \sum e^1 r \alpha(e^2) = \sum re^1 \alpha(e^2)
\} = Z, $$
since $\supseteq$ is clear, and $\subseteq$ follows from taking $\alpha = \id_A$ and
observing $\sum e^1 r e^2 \in Z$.  
The Morita context bimodule
$$ Q = \{ q \in \cC \, |\, \forall  \alpha \in \cC, \,  
\sum q(-e^1) \alpha(e^2) = \sum q(-) e^1 \alpha(e^2) \}
= \lambda(R), $$
since $\supseteq$ is clear, and $\subseteq$ follows from taking $\alpha = \id_A$,
whence $q = \sum q(-e^1)e^2 =\sum  \lambda(q(e^1)e^2) \in \lambda(R)$. 

That $(Z,\cC, Q,R,\sigma, \tau)$ is a Morita context is now
straightforward; $\tau$ being epi by an old lemma of Hirata and Sugano \cite{HS}.

Recall that $B \rightarrow A$ is H-separable (after Hirata) if there
are (Casimir) elements $e_i \in (A \ot_B A)^A$ and $r_i \in R$ (the centraliser)
such that $1_A \ut 1_A = \sum_i e_i r_i (= \sum_i r_ie_i)$ (a very strong version
of Eqs.~(\ref{eqs:depthtwo}) above). It is well-known that $A$ is a separable extension of
$B$.  Moreover, 
$$R \ot_Z R^{\rm op} \stackrel{\cong}{\longrightarrow} \lrend BBA, $$
via $r \ut r' \mapsto \lambda(r) \rho(r')$ with inverse
$\alpha \mapsto \sum_i \alpha(e^1_i)e^2_i \ut r_i$.
Whence $\sigma$ is an isomorphism if we begin with an H-separable extension.  
\end{proof}

If $A$ is a separable algebra over a (commutative) ground ring $B$, then the proposition shows
that the center $Z$ and $\rend B A$ (with the exotic multiplication above)
are related by a Morita context, which is strict if $Z = B 1_A$,
i.e.\ $A$ is Azumaya. 

As a third example, we may instead work with the dual bialgebroid in \cite{KadSzlach:depth}
and prove that a split D2 extension $B \rightarrow A$ has a coseparable
$R$-coring structure on $(A \ot_B A)^B$ which is essentially a restriction
of $\cC$ in Example~\ref{ex:split}. 

\section{Are biseparable corings  Frobenius?}

In this section we will show that a one-sided, slightly stronger
 version of the problem in \cite{CaeKad:bis} is equivalent
to the problem if cosplit, coseparable corings with a condition of finite
projectivity are Frobenius.  Given the techniques developed for corings and
the many examples coming from entwined structures \cite{Brz:str}, we expect
this equivalence to be useful in solving this problem. 

As recalled in Section~2, an $A$-coring $\cC$ is coseparable if the forgetful functor
$F:  \M^{\cC} \to \M_A$ is separable (cf.\ \cite{NasBer:sep} for the 
definition of a separable functor). Dually, we say that $\cC$ is \textit{cosplit} if
the functor $- \otimes_A \cC$ is a separable functor  from the category
of right $A$-modules $\M_A$ into the category of right $\cC$-comodules
$\M^{\cC}$. (Recall that $F$ is the left adjoint of $-\otimes_A \cC$ \cite{Guz:coi}.)

An $A$-coring $\cC$ determines two ring extensions $\iota^*: A \to 
\cC^*$ and ${}^*\iota: A \to {}^*\cC$ where $\cC^* := \rhom A \cC A$ and ${}^*\cC :=
\lhom A \cC A$, i.e., the right and left duals of $\cC$. The ring structure
on $\cC^*$ is given by $(\xi\xi')(c) = \sum \xi(\xi'(c \sw 1) c \sw 2)$
($\xi,\xi' \in \cC^*, c \in \cC$) with unity $\eC$
and the natural Abelian group structure, 
while the ring structure on ${}^*\cC$ is given by $(\xi\xi')(c) = \sum 
\xi'(c \sw 1 \xi(c \sw 2))$.
The mappings $\iota^*$ and ${}^*\iota$ are given by $\iota^*(a) = \eC(a-) (= a\eC)$
and ${}^*\iota(a) = \eC(-a)$.
  We note by short calculations that
the induced $(A,A)$-bimodule structures on $\cC^*$ and ${}^*\cC$ coincide with the usual
structures, which we recall are given by $(a\xi a')(c) := a\xi (a'c)$ for 
$\xi\in \cC^*$ and $a, a' \in A$ and
$(a\xi a')(c) = \xi (ca)a'$ for $\xi \in {}^*\cC$. Also note that
${}^*\iota$ and $\iota^*$ are monomorphisms if $\eC: \cC \to A$ is surjective.

Recall from \cite{CaeKad:bis} that a ring extension $B \to A$ is biseparable if
it is split, separable and the natural modules $A_B$ and ${}_BA$ are finitely generated  projective.
We will say that $B \to A$ is \textit{left or right biseparable} if $B \to A$ is split, separable but 
only one of  ${}_BA$
or $A_B$, respectively, need be finitely generated  projective.
This motivates the following 
\begin{definition}
An $A$-coring $\cC$ is said to be {\em biseparable} if $\cC_A$ and ${}_A\cC$ are finitely
generated projective and $\cC$ is cosplit as well as coseparable.
\end{definition}

\begin{proposition}
\label{prop-canon}
If $B \to A$ is a biseparable extension, then the canonical Sweedler $A$-coring $\cC := A \otimes_B A$
is a biseparable coring.
\end{proposition}
\begin{proof}
Since $B \to A$ is separable, the induction functor $- \otimes_A \cC$ from
$\M_A$ into $\M^{\cC}$ is separable by \cite[Corollary~3.4]{Brz:str}.
Since $B \to A$ is split, and $A_B$ is a projective
generator (therefore faithfully flat), the forgetful functor $F: \M^{\cC} \to
\M_A$ is a separable functor by \cite[Corollary~3.7]{Brz:str}. It follows by definition that the canonical
coring $\cC$ is cosplit and coseparable.

Finally we note that ${}_BA$ finitely generated  projective implies ${}_AA \ut_B A$ finitely generated  projective.
Similarly, $\cC_A$ is finitely generated  projective, and we conclude that $\cC$ is biseparable. 
\end{proof} 

Recall that an $A$-coring $\cC$ is said to be {\em Frobenius} if the forgetful functor
$F: \M^{\cC} \to \M_A$ is a Frobenius functor (has the same left and right adjoint),
i.e. $-\otimes_A \cC$ is also a left adjoint of $F$ \cite{Brz:str,Brz:tow}. 
Motivated by the question in \cite{CaeKad:bis} let us make

\begin{conjecture}
\label{con-coring}  
A biseparable $A$-coring $\cC$ is Frobenius.
\end{conjecture}

\begin{proposition}
\label{prop-two}
If Conjecture~\ref{con-coring} is true, then biseparable extensions are Frobenius.
\end{proposition}
\begin{proof}
Given a biseparable extension $B \to A$, its canonical coring $\cC  = A 
\ot_BA$ is
biseparable by the previous proposition.  If $\cC$ is then a Frobenius coring by hypothesis,
it follows from \cite[Theorem~2.7]{Brz:tow} that $B \to A$ is a Frobenius extension, since
${}_BA$ is faithfully flat.
\end{proof}

We now proceed to establish a converse to this proposition.  

\begin{proposition}
\label{prop-cosplit}
If $\cC$ is a cosplit $A$-coring, then $\iota^*: A \to \cC^*$ and ${}^*\iota: A \to
{}^*\cC$ are both split extensions.
\end{proposition}
\begin{proof}
By \cite[Theorem 3.3]{Brz:str}, $\cC$ is cosplit 
iff there is $e \in \cC^A$ such that $\eC(e) = 1$. (In other words,
$\eC: \cC \to A$ is a split $(A,A)$-epimorphism.)
  We now define a ``conditional expectation'' or bimodule projection $E^*: \cC^* \to
A$, respectively ${}^*E: {}^*\cC \to A$ simply by
\[
E^*(\xi) = \xi (e), \ \ \ {}^*E(\xi') = \xi'(e) \ \ \ (\xi\in \cC^*, 
\xi' \in {}^*\cC)
.
\]
Note that $E^*(\eC) = 1_A = {}^*E(1_{{}^*\cC})$ and
$$ E^*(a\xi a') = a\xi (a'e) = a\xi (ea') = a\xi (e)a' = aE^*(\xi)a'$$
for $a,a' \in A$, whence $E^*$ and similarly ${}^*E$ give $(A,A)$-bimodule
splittings of $\iota^*$ and ${}^*\iota$. 
\end{proof}
 
For example, the canonical Sweedler coring $\cC$ of a ring extension $B \to A$ is
cosplit if and only if $B \to A$ is a separable extension. Now $\cC^* \cong \rend B A$ as rings 
via $\xi \mapsto \xi(- \ut 1)$ with inverse $$ f \longmapsto [a \ut a' \mapsto 
f(a)a'].$$
Since $\iota^*$ corresponds to the left regular representation $\lambda: A
\to \rend B A$, we recover  results by M\"uller and Sugano
 that $\lambda$ is a split extension if $B \to A$ is separable.

\begin{proposition}
\label{prop-coseparable}
Let $\cC$ be an $A$-coring. If ${}_A\cC$ is finitely generated  projective,
then  $\cC$ is coseparable if and only if 
${}^*\iota: A \to {}^*\cC$ is a separable extension.
If $\cC_A$ is finitely generated  projective, then 
$\cC$ is a coseparable $A$-coring if and only if
 $\iota^*: A \to \cC^*$ is a separable extension.  
\end{proposition}
\begin{proof}
We will prove the first statement, the second follows similarly.
The category of right comodules $\M^{\cC}$ is isomorphic to the
category $\M_{{}^*\cC}$ of right modules over ${}^*\cC$ \cite[Lemma~4.3]{Brz:str}.
Recall that
 given a coaction $\varrho^{M}: M_A \to M \ut_A \cC$, we 
define an action of $\xi \in {}^*\cC$ on $m \in M \in \M^{\cC}$
 by $m \cdot \xi = \sum m \sw 0 \xi (m \sw 1)$. It is trivial to check that
$(M, \cdot) \in \M_{{}^*\cC}$.
Inversely, given dual bases $\{ \xi_i  \in {}^*\cC\}$ and $\{ c_i \in \cC \}$
such that $c = \sum_i \xi_i(c) c_i$ for each $c \in \cC$, 
and right action of ${}^*\cC$ on $M \in \M_{{}^*\cC}$, we define
a coaction
$$ \varrho^{M}(m ) = \sum_i m \cdot \xi_i \otimes_A c_i. $$
It is easily checked that $(M, \varrho^{M}) \in \M^{\cC}$,
and that the two operations are natural and inverses to one another,
so that $\M_{{}^*\cC} \cong \M^{\cC}$. 

Now $\cC$ is coseparable if and only if the forgetful functor $F: \M^{\cC} \to
\M_A$ is a separable functor. Since
$$ m \cdot {}^*\iota(a) = \sum m \sw 0 \eC(m \sw 1)a = ma $$
for each $a \in A, m \in M \in \M^{\cC}$, 
$F$ corresponds
under the isomorphism of categories above to the forgetful functor $G:
\M_{{}^*\cC} \to \M_A$ induced by ${}^*\iota: A \to {}^*\cC$.  
 But $G$ is a separable functor if and only if ${}^*\iota$ is a 
separable extension \cite{NasBer:sep}.
\end{proof}

\begin{proposition}
\label{prop-frobenius}
Suppose $\cC$ is an $A$-coring which is  reflexive as a left and right  $A$-module. Then
$\cC$ is a Frobenius coring  if and only if 
 ${}^*\iota: A \to {}^*\cC$ is a Frobenius extension
 if and only if 
$\iota^*: A \to \cC^*$  is a Frobenius extension.
\end{proposition}
\begin{proof}
The proof is quite similar to the proof of Proposition~\ref{prop-coseparable}
(cf. \cite[Theorem~4.1]{Brz:str}). If ${}^*\iota$ is Frobenius, it
follows that ${}^*\cC_A$ is finitely generated  projective, so ${}_A({}^*\cC)^*$ is finitely generated  projective,
whence by reflexivity ${}_A\cC$ is finitely generated  projective.  Then
 the categories $\M_{{}^*\cC}$
and $\M^{\cC}$ are isomorphic.  But the forgetful functor $G:
\M_{{}^*\cC}\to \M_A$ has equal left and right adjoint if
 ${}^*\iota$ is Frobenius, in which case $F$ is Frobenius and $\cC$ is
a Frobenius coring. The other case is entirely similar. 
Conversely, if $\cC$ is a Frobenius $A$-coring, then both $\cC_A$ and
${}_A\cC$ are finitely generated  projective \cite[Corollary~2.3]{Brz:tow}. The rest follows from the functorial
definitions of Frobenius extension and coring applied to either isomorphism of left or right
module and comodule categories .    
\end{proof}

\begin{theorem}
\label{thm-equivalence}
Biseparable corings are Frobenius if and only if left or right biseparable extensions are Frobenius.
\end{theorem}
\begin{proof}
Without one-sidedness, we saw $\Rightarrow$ in Proposition~\ref{prop-two}.  ($\Leftarrow$)  Suppose
$\cC$ is  a biseparable $A$-coring.  Then ${}^*\iota: A \to {}^*\cC$ and
$\iota^*: A \to \cC^*$ are split, separable extensions by Propositions~\ref{prop-cosplit}
and~\ref{prop-coseparable}. Since ${}_A\cC$ and $\cC_A$ are finitely generated  projective,
it follows that ${}^*\cC_A$ and ${}_A\cC^*$ are finitely generated  projective. If either
left or right biseparable extensions are Frobenius, then either 
$\iota^*$ or ${}^*\iota$ is a Frobenius extension.  In either case,
Proposition~\ref{prop-frobenius} shows $\cC$ to be a Frobenius coring.
\end{proof}     

We note the following special ``depth one'' case for which
there is a  solution to our conjecture.
If $\cC$ is a centrally projective $A$-bimodule, i.e.,
as $A$-$A$-bimodules $\cC \oplus W \cong \oplus^n A$ for some $(A,A)$-bimodule $W$,
and $\cC$ is moreover biseparable, then $\cC$ is Frobenius by
 a classical result of Sugano:
\begin{proposition}
If $\cC$ is a centrally projective, cosplit, coseparable $A$-coring,
then $\cC$ is Frobenius.
\end{proposition}
\begin{proof}
We easily obtain $\cC^* \oplus W^* \cong \oplus^n A$ as $A$-$A$-bimodules,
whence
 $\iota^*: A \to \cC^*$
is a centrally projective separable extension by Proposition~\ref{prop-coseparable},
and monomorphism since $\eC$ is surjective.
By 
\cite[Theorem 2]{Sug:sep} $\iota^*$ is a Frobenius extension.
Then by Proposition~\ref{prop-frobenius},
 $\cC$ is a Frobenius coring.  
\end{proof}

\section*{Acknowledgements}
Tomasz Brzezi\'nski  would like to thank the Engineering and Physical Sciences Research 
Council for an Advanced Fellowship. Lars Kadison thanks Dmitri Nikshych and the University
of New Hampshire Department of Mathematics and Statistics for their generosity.


\begin{thebibliography}{Bibliography}{}
    \bibitem{Abu:Mor} J.Y.\ Abuhlail. Morita contexts for corings 
    and equivalences. Preprint, 2002.
    \bibitem{Bas:KTh} H.\ Bass {\em Algebraic K-Theory}. W.A.\ 
    Benjamin, inc., New York, 1968.
    \bibitem{BerHau:cog} G.M.\ Bergman and A.O.\ Hausknecht {\em 
    Cogroups and Co-rings in Categories of Associative Rings}. AMS, 
    Providence R.I., 1996.
    \bibitem{BerKea:cat} A.J.\ Berrick and M.A.\ Keating. {\em 
    Categories and Modules with K-Theory in View.} Cambridge 
    University Press, 2000.

\bibitem{Brz:str} T. Brzezi\'nski. The structure of corings. Induction functors, 
Maschke-type theorem, and Frobenius and 
    	Galois-type properties. {\em Alg.\ 
    Rep.\ Theory} in press.  Preprint arXiv: math.RA/0002105, 2000.
\bibitem{Brz:tow} T. Brzezi\'nski. Tower of corings. {\em Commun.\ 
    Algebra}, in press. Preprint arXiv: math.RA/0201014, 2002.
\bibitem{BrzMil:bia} T.\ Brzezi\'nski and G.\ Militaru. Bialgebroids, $\times_A$-bialgebras and duality. {\em J.\ Algebra}, 251: 279--294, 2002.
\bibitem{CaeKad:bis} S. Caenepeel and L. Kadison. Are biseparable extensions Frobenius? 
{\em   K-Theory}, 24: 361--383, 2001.
 \bibitem{Cae:Mor} S.\ Caenepeel, J.\ Vercruysse and S.\ Wang. Morita 
 theory for corings and cleft entwining structures. Preprint arXiv: 
 math.RA/0206198, 2002.
 \bibitem{GomLou:cos} J.\ G\'omez-Torrecillas and A.\ Louly. 
 Coseparable corings.  Preprint arXiv: math.RA/0206175, 2002.

 \bibitem{Guz:coi} F. Guzman. Cointegrations, relative cohomology for 
  comodules and coseparable corings. {\em J.\ Algebra}, 126:211--224, 1989.

\bibitem{HS} K.\ Hirata and K.\ Sugano, 
On semisimple extensions and separable extensions over non commutative
rings, {\em J. Math.\ Soc.\ Japan}, 18: 360--373, 1966.

\bibitem{KadSzlach:depth} L. Kadison and K. Szlach\'anyi,
Bialgebroid  actions on depth two extensions and duality. 
           Preprint arXiv: math.RA/0108067, 2001.
\bibitem{Lu:Hop}
J.H. Lu. Hopf algebroids and quantum groupoids,
{\em Int. J. Math.}, 7:47--70, 1996.

\bibitem{NasBer:sep} C. Nastasescu, M. Van den Bergh and F. Van Oystaeyen. Separable functors applied to graded rings. {\em J.\ 
    Algebra}, 123:397--413, 1989.
\bibitem{Q:notes} D. Quillen. Module theory over non-unital rings. {\em Preprint} 1997. 
\bibitem{Pie:alg} R.S. Pierce. Associative Algebras. {\em Springer Verlag} 1982.
\bibitem{Raf:sep} M.D. Rafael. Separable functors revisited. {\em Commun.\ 
    Algebra}, 18:1445--1459, 1990.
\bibitem{Sug:sep} K. Sugano. Separable extensions and Frobenius extensions. {\em Osaka J.  
    Math.}, 7:291--299, 1970.
\bibitem{Swe:pre}
M.\ Sweedler, The predual theorem to the Jacobson-Bourbaki theorem,
{\em Trans.\ Amer.\ Math.\ Soc.} 213:391--406, 1975.

    \bibitem{Tay:big} J.L.\ Taylor. A bigger Brauer group. {\em Pac.\ 
    J.\ Math.}, 103:163--203, 1962.
    
\end{thebibliography}
\end{document}